\documentclass[twocolumn]{autart}    

\usepackage{graphicx}          
                               
\usepackage{amsmath}
\usepackage{amssymb}

\usepackage{graphics} 
\usepackage{epsfig} 
\usepackage{epstopdf}
\usepackage{subfigure}
\usepackage{color}

\graphicspath{{./Figures/}}
\allowdisplaybreaks[4]

\def\R{\mathbb{R}}
\def\N{\mathbb{N}}

\def\C{\mathbb{C}} 

\usepackage{bbm}

\begin{document}

\begin{frontmatter}

\title{Integral action for setpoint regulation control of a reaction-diffusion equation in the presence of a state delay\thanksref{footnoteinfo}} 

\thanks[footnoteinfo]{This publication was supported in part by a research grant from Science Foundation Ireland (SFI) under grant number 16/RC/3872 and is co-funded under the European Regional Development Fund and by I-Form industry partners. 
The work of the first author was supported by ANR PIA funding: ANR-20-IDEES-0002.
Corresponding author H.~Lhachemi.}

\author[CS]{Hugo Lhachemi}\ead{hugo.lhachemi@centralesupelec.fr}, 
\author[UCD]{Ammar Malik}\ead{ammar.malik@ucdconnect.ie}, 
\author[UCD,ICL]{Robert Shorten}\ead{r.shorten@imperial.ac.uk},               

\address[CS]{Universit{\'e} Paris-Saclay, CNRS, CentraleSup{\'e}lec, Laboratoire des signaux et syst{\`e}mes, 91190, Gif-sur-Yvette, France}  
\address[UCD]{School of Electrical and Electronic Engineering, University College Dublin, Dublin, Ireland}  
\address[ICL]{Dyson School of Design Engineering, Imperial College London, London, U.K}             


\begin{abstract}                          
This paper is concerned with the regulation control of a one-dimensional reaction-diffusion equation in the presence of a state-delay in the reaction term. The objective is to achieve the PI regulation of the right Dirichlet trace with a command selected as the left Dirichlet trace. The control design strategy consists of the design of a PI controller on a finite dimensional truncated model obtained by spectral reduction. By an adequate selection of the number of modes of the original infinite-dimensional system, we show that the proposed control design procedure achieves both the exponential stabilization of the original infinite-dimensional system as well as the setpoint regulation of the right Dirichlet trace.  
\end{abstract}

\end{frontmatter}

\section{Introduction}

The proportional integral (PI) regulation control of infinite-dimensional systems, and in particular of partial differential equations (PDEs), has attracted much attention in the recent years. Early works dealt with bounded control operators~\cite{pohjolainen1982robust} while the extension to the case of unbounded control operators was reported in~\cite{xu1995robust}. The last decade has seen an intensification of the efforts in this research direction. PI boundary control of linear hyperbolic systems~\cite{bastin2015stability,dos2008boundary,lamare2015control,xu2014multivariable}, as well as the extension to nonlinear transport equations~\cite{coron2019pi,trinh2017design} have been reported. Other types of PDEs have also been studied. This includes reaction-diffusion equations~\cite{lhachemi2019pi}, wave equations used to model drilling systems~\cite{barreau2019practical,terrand2018regulation}, as well as semilinear wave equations~\cite{lhachemi2020pi}. The possible addition of an integral action to open-loop stable semigroups was investigated in~\cite{terrand2019adding}. 

We study in this paper the boundary PI regulation control of a reaction-diffusion equation in the presence of a state-delay in the reaction term. Since delays are ubiquitous in practical applications, the topic of boundary stabilization of PDEs in the presence of delays, either in the control input~\cite{krstic2009control,lhachemi2019feedback,lhachemi2019lmi,nicaise2008stabilization,prieur2018feedback} or in the state~\cite{hashimoto2016stabilization,kang2017boundary,lhachemi2020boundary}, has also attracted much attention in the recent years. However, it is worth noting that none of the aforementioned works embracing PI control design for PDEs was concerned with the possible presence of state-delays. This paper is a first step in that research direction. Specifically, the objective of this work is to extend the result reported in~\cite{lhachemi2020boundary}, which solely dealt with the boundary stabilization of a reaction-diffusion equation in the presence of a state-delay, to the PI regulation control of a Dirichlet trace. More precisely we consider the PDE: 
\begin{subequations}\label{eq: heat equation - PDE}
\begin{align}
& y_t(t,x) = a y_{xx}(t,x) + b y(t,x) + c y(t-h(t),x) , \\
& y(t,0) = u(t) \label{eq: heat equation - PDE - control input} , \\
& \cos(\theta) y(t,1) + \sin(\theta) y_x(t,1) = 0 , \label{eq: heat equation - PDE - BC RHS} \\
& y(\tau,x) = \phi(\tau,x) , \quad \tau\in[-h_M,0]
\end{align}
\end{subequations}
for $t > 0$ and $x \in (0,1)$. Here $a>0$, $b,c \in\R$ with $c \neq 0$,  $h : \R_+ \rightarrow [h_m , h_M]$ is continuous with $0 < h_m < h_M$, and $\theta \in (0,\pi/2)$. The state at time $t$ is $y(t,\cdot):[0,1]\rightarrow\R$. The control input is $u(t) \in \R$ and applies to the left Dirichlet trace (\ref{eq: heat equation - PDE - control input}). On the right-hand side of the domain, we consider the Robin boundary condition (\ref{eq: heat equation - PDE - BC RHS}). The initial condition is $\phi : [-h_M,0] \times (0,1) \rightarrow \R$. The control objective is to design a PI controller in order to stabilize (\ref{eq: heat equation - PDE}) while achieving the setpoint regulation control of the right Dirichlet trace $z(t) = y(t,1)$. In particular, denoting by $r: \R_+ \rightarrow \R$ a continuous reference signal, $y(t,1)$ must achieve the setpoint tracking of $r(t)$.

The strategy for solving the above control design problem goes as follows. Inspired by~\cite{coron2004global}, a finite dimensional-truncated model is obtained by spectral reduction. The order of the state-delayed truncated model is selected to ensure the stability of the residual infinite-dimensional dynamics. Then, inspired by~\cite{lhachemi2019pi} but with the challenge of a state-delayed term, the truncated model is augmented with an integral component to ensure the setpoint tracking of $z(t)$. Finally, the feedback law is obtained by pole shifting. We assess the exponential stability of the closed-loop system, as well as the setpoint regulation control of the right Dirichlet trace. In the presence of an additive boundary perturbation in the control input, we show that the closed-loop system is exponentially input-to-state stable (ISS). This objective requires to work simultaneously with the original representation of the plant (for ISS purposes w.r.t. boundary disturbances) and an homogeneous version of the PDE (to analyze the system output) while handling the state-delay for both stability and setpoint regulation assessment.

The control design strategy is introduced in Section~\ref{sec: control design strategy}. The equilibrium conditions of the closed-loop system and the related dynamics of deviations are presented in Section~\ref{sec: eq cond and dynamics of deviations}. The stability analysis is reported in Section~\ref{sec: stability analysis} while the reference tracking assessment is completed in Section~\ref{sec: setpoint regulation assessment}. The robustness of the control strategy w.r.t. delay mismatches is studied in Section~\ref{sec: robustness}. Finally, numerical simulations are carried out in Section~\ref{sec: simulation} while concluding remarks are formulated in Section~\ref{sec: conclusion}.

\section{Control design strategy}\label{sec: control design strategy}

\subsection{Spectral reduction and truncated model}

Let $\mathcal{H} = L^2(0,1)$ with the inner product $\left< f , g \right> = \int_0^1 fg \,\mathrm{d}x$. System (\ref{eq: heat equation - PDE}) can be rewritten as 
\begin{subequations}\label{eq: heat equation - ACS}
\begin{align}
& \dfrac{\mathrm{d}X}{\mathrm{d}t}(t) = \mathcal{A}X(t) + c X(t-h(t)) , \label{eq: heat equation - ACS - ODE} \\
& \mathcal{B}X(t) = u(t) , \\
& X(\tau) = \Phi(\tau) = \phi(\tau,\cdot)  , \quad \tau\in[-h_M,0] 
\end{align}
\end{subequations}
for $t \geq 0$ with $\mathcal{A} : D(\mathcal{A}) \subset \mathcal{H} \rightarrow \mathcal{H}$ defined on $D(\mathcal{A}) = \left\{ f \in H^2(0,1) \,:\, \cos(\theta) f(1) + \sin(\theta) f'(1) = 0  \right\}$, with $\theta \in (0,\pi/2)$,  by $\mathcal{A} f = a f'' + b f$ and the boundary operator $\mathcal{B} : D(\mathcal{B}) \subset \mathcal{H} \rightarrow \R$ defined on $D(\mathcal{B}) = H^1(0,1)$ by $\mathcal{B} f = f(0)$. We define the disturbance free operator $\mathcal{A}_0 = \left.\mathcal{A}\right\vert_{D(\mathcal{A}_0)}$ on $D(\mathcal{A}_0) = D(\mathcal{A}) \cap \operatorname{ker}(\mathcal{B})$. It is well-known that $\mathcal{A}_0$ generates a $C_0$-semigroup. We introduce $L \in \mathcal{L}(\R,\mathcal{H})$ defined for any $u \in \R$ by $[Lu](x) = (1-x)^2 u$, $x \in [0,1]$. $L$ has been selected such that its range satisfies $\operatorname{R}(L) \subset D(\mathcal{A})$ and $\mathcal{B}L = I_\R$. Hence, following the terminology of~\cite[Sec.~3.3]{curtain2012introduction}, the pair $(\mathcal{A},\mathcal{B})$ defines a boundary control system with associated lifting operator $L$. We define $\mathcal{A}_{c} \triangleq \mathcal{A} + c I_\mathcal{H}$ and $\mathcal{A}_{c,0} \triangleq \mathcal{A}_0 + c I_\mathcal{H}$ on $D(\mathcal{A}_{c}) = D(\mathcal{A})$ and $D(\mathcal{A}_{c,0}) = D(\mathcal{A}_0)$, respectively. From the Sturm-Liouville theory, it well known that the eigenvalues of $\mathcal{A}_{c,0}$ are simple and form a decreasing sequence $(\lambda_n)_{n \geq 0} \in\R^\N$ with $\lambda_n \rightarrow - \infty$ when $n \rightarrow + \infty$. Moreover, one can select the associated eigenvectors such that $(e_n)_{n \geq 0}$ forms a Hilbert basis of $\mathcal{H}$. Using the terminology of~\cite[Def.~2.3.4]{curtain2012introduction}, $\mathcal{A}_{c,0}$ is a Riesz spectral operator: $D(\mathcal{A}_{c,0}) = \left\{ f \in L^2(0,1) \,:\, \sum_{n\geq0} \vert\lambda_n\vert^2 \vert \left< f , e_n \right> \vert^2 < \infty \right\}$ and $\mathcal{A}_{c,0} f = \sum_{n \geq 0} \lambda_n \left< f , e_n \right> e_n$ for all $f \in D(\mathcal{A}_{c,0})$. Standard computations give $\lambda_n = b + c - a r_n^2$ and $e_n = 2 \sqrt{\frac{r_n}{2 r_n - \sin(2 r_n)}} \sin(r_n \cdot)$ with $n \in \N$ where $r_n > 0$ is the unique number $r \in (n\pi , (n+1)\pi)$ such that $r \cot(r) = - \cot(\theta)$. This yields $\lambda_n \sim -an^2\pi^2$ and $e_n(1) = O(1)$ as $n \rightarrow +\infty$. Introducing $x_n(t) = \left< X(t) , e_n \right>$ the coefficients of projection of the system trajectory into the Hilbert basis, we have that $X(t) = \sum_{n\geq0} x_n(t) e_n$ and $\Vert X(t) \Vert^2 = \sum_{n\geq0} \vert x_n(t) \vert^2$. Assuming that\footnote{Such a regularity for the forthcoming control law will be assessed in the sequel.} $u$ is continuously differentiable and $\Phi$ is continuous, the mild solution $X \in \mathcal{C}^0(\R_+;\mathcal{H})$ of (\ref{eq: heat equation - ACS}) is such that $x_n$ is continuously differentiable and (see~\cite{lhachemi2020boundary} for details)
\begin{equation}\label{eq: ODE coefficient of projection}
\dot{x}_n(t) = \lambda_n x_n(t) + c \{ x_n(t-h(t)) - x_n(t) \} + (a_n + \lambda_n b_n) u(t) 
\end{equation}
with 
\begin{equation}\label{eq: def an bn}
a_n = \left< \mathcal{A}_{c} L \mathbbm{1} , e_n \right> , \quad
b_n = - \left< L \mathbbm{1} , e_n \right> 
\end{equation}
where $\mathbbm{1}$ denotes here the unit element of $\R$. Note that due to the presence of the state-delay, there may exist delays for which certain modes $x_n$, hence the PDE, are unstable even if $\lambda_n - c = b - a r_n^2 < 0$ and $c < 0$ provided $c$ is large enough~\cite[Sec.~3.3]{michiels2002continuous}. For a given integer $N \geq 0$ selected such that $\lambda_{n} < 0$ for all $n \geq N+1$ and which will be further constrained later, we define the followings:
\begin{subequations}\label{eq: components of truncate model}
\begin{align}
Y(t) &= \begin{bmatrix} x_0(t) & \ldots & x_{N}(t) \end{bmatrix}^\top \in \R^{N+1} , \\
Y_\Phi(\tau) & = \begin{bmatrix} \left< \Phi(\tau) , e_0 \right> & \ldots & \left< \Phi(\tau) , e_N \right> \end{bmatrix}^\top \in \R^{N+1} , \\
A & = \mathrm{diag}(\lambda_n)_{0 \leq n \leq N} \in \R^{(N+1) \times (N+1)} , \\
B & = (a_n + \lambda_n b_n)_{0 \leq n \leq N} \in \R^{(N+1)} .
\end{align}
\end{subequations}
Then we obtain the truncated model:
\begin{subequations}\label{eq: truncate model}
\begin{align}
\dot{Y}(t) & = A Y(t) + c \{ Y(t-h(t)) - Y(t) \} + Bu(t) \\
Y(\tau) & = Y_\Phi(\tau) , \quad \tau\in[-h_M,0]
\end{align}
\end{subequations}

\subsection{Addition of an integral component}

The objective is now to augment the truncated model (\ref{eq: truncate model}) with an integral component to achieve the setpoint regulation control of the right Dirichlet trace $z(t) = y(t,1)$. 

We first need to express the right Dirichlet trace $y(t,1)$ in function of the coefficients of projection $x_n$.

\begin{lem}\label{lem: seies expansion of f(1)}
Let $\theta \in (0,\pi/2)$. For all $f \in D(\mathcal{A}_{c,0})$ we have $f(1) = \sum_{n \geq 0} \left< f , e_n \right> e_n(1)$.
\end{lem}

The proof of this Lemma, which is omitted, essentially relies on the Riesz-spectral property of $\mathcal{A}_{c,0}$. We cannot directly apply the above series expansion to the trajectory $X$ of our system because, in general, $X(t) \notin D(\mathcal{A}_{c,0})$. However, if we assume that $X \in \mathcal{C}^0(\R_+;D(\mathcal{A})) \cap \mathcal{C}^1(\R_+;\mathcal{H})$ is a classical solution of (\ref{eq: heat equation - ACS}), one has $W(t) = X(t) - L u(t) \in D(\mathcal{A}_{0}) = D(\mathcal{A}_{c,0})$ with in particular $y(t,1) = [X(t)](1) = [W(t)](1)$. Hence, introducing $w_n(t) = \left< W(t) , e_n \right> = x_n(t) + b_n u(t)$, we obtain that $y(t,1) = [X(t)](1) = \sum_{n \geq 0} w_n(t) e_n(1)$ for all $t \geq 0$. Since $u$ is of class $\mathcal{C}^1$ (we will actually need $u$ of class $\mathcal{C}^2$ to ensure the existence of classical solutions), we have that $w_n$ is of class $\mathcal{C}^1$ and, from (\ref{eq: ODE coefficient of projection}), 
\begin{align}
\dot{w}_n(t) 
& = \lambda_n w_n(t) + c \{ w_n(t-h(t)) - w_n(t) \} + a_n u(t) \nonumber \\
& \phantom{=}\; - c b_n \{ u(t-h(t)) - u(t) \} + b_n \dot{u}(t) \label{eq: spectral reduction EDO wn}
\end{align}
for $t \geq h_M$.

Remind that our objective is to achieve the setpoint regulation control of the system output $z(t) = y(t,1)$. In order to introduce in a comprehensive manner the proposed integral component $\zeta(t) \in\R$ that will be used to augment the truncated model (\ref{eq: truncate model}), consider first the classical integral component given by $\dot{\chi}(t) = y(t,1) - r(t) = \sum_{n \geq 0} w_n(t) e_n(1) - r(t)$ for $t \geq 0$. Here $r(t)\in\R$ stands for a reference signal. Recall that the second equality holds only when considering classical solutions for (\ref{eq: heat equation - ACS}). As the above series expansion involves all the modes of the system, and in particular the coefficients of projection $w_n(t)$ for $n \geq N+1$, the integral component $\chi$ cannot be directly included into the dynamics of the truncated model (\ref{eq: truncate model}). To solve this issue, we introduce the following preliminary change of variable $\zeta_p(t) = \chi(t) + \sum_{n \geq N+1} \frac{e_n(1)}{\lambda_n} \{ b_n u(t) - w_n(t) \}$. Note that the convergence of the series follow from $\lambda_n \sim -an^2\pi^2$ and $e_n(1) = O(1)$ as $n \rightarrow +\infty$. Based on (\ref{eq: spectral reduction EDO wn}) we obtain, for $t \geq h_M$, $\dot{\zeta}_p(t) = \sum_{n = 0}^{N} x_n(t) e_n(1) + \alpha u(t) - r(t) - c \sum_{n \geq N+1} \frac{e_n(1)}{\lambda_n} \{ w_n(t-h(t)) - w_n(t) \} + c \sum_{n \geq N+1} \frac{e_n(1)}{\lambda_n} b_n \{ u(t-h(t)) - u(t) \}$
where
\begin{equation}\label{eq: def alpha}
\alpha = \sum\limits_{n = 0}^{N} b_n e_n(1) - \sum\limits_{n \geq N+1} \dfrac{a_n}{\lambda_n} e_n(1) .
\end{equation}
We now note that the two last terms of the above identity describing the $\zeta_p$-dynamics have a null contribution at equilibrium. This observation motivates the introduction of the below $\zeta$-dynamics. Assuming that the delay $h$ is known (robustness w.r.t. delay mismatches will be discussed later in Section~\ref{sec: robustness}), we mimic the structure of the dynamics of the truncated model (\ref{eq: truncate model}) by defining for $t \geq 0$ the integral component $\zeta(t) \in\R$ as follows:
\begin{subequations}\label{eq: zeta dynamics}
\begin{align}
\dot{\zeta}(t) & = \sum\limits_{n = 0}^N x_n(t) e_n(1) + c \{ \zeta(t-h(t)) - \zeta(t) \} \\
& \phantom{=}\; + \alpha u(t) - r(t) , \nonumber \\
\zeta(\tau) & = \zeta_0(\tau) , \quad \tau \in [-h_M,0]
\end{align}
\end{subequations}

\begin{rem}
The $\zeta$-dynamics achieves the same equilibrium condition as the $\zeta_p$-dynamics. As we will show later in Section~\ref{sec: eq cond and dynamics of deviations}, the integral component (\ref{eq: zeta dynamics}) ensures that the equilibirum condition $(X_e,\zeta_e)$ of the forthcoming closed-loop system, associated with some constant reference signal $r(t) = r_e$, achieves the desired reference tracking for the right Dirichlet trace, i.e., $X_e(1) = r_e$.
\end{rem}

\begin{rem}
Even if (\ref{eq: zeta dynamics}) has been motivated and derived by considering classical solutions of (\ref{eq: heat equation - ACS}), the dynamics (\ref{eq: zeta dynamics}) actually makes sense for any mild solutions of (\ref{eq: heat equation - ACS}). 
\end{rem}

Since (\ref{eq: zeta dynamics}) only involves the $N+1$ first modes of the system, we can now augment the dynamics of the truncated model (\ref{eq: truncate model}) with the $\zeta$-dynamics as follows:
\begin{subequations}\label{eq: aug truncate model}
\begin{align}
\dot{Y}_a(t) & = A_a Y_a(t) + c \{ Y_a(t-h(t)) - Y_a(t) \} \\
& \phantom{=}\; + B_a u(t) + \Gamma(t) , \nonumber \\
Y_a(\tau) & = Y_{\Phi,a}(\tau) , \quad \tau\in[-h_M,0]
\end{align}
\end{subequations}
where $C = \begin{bmatrix} e_0(1) & \ldots & e_N(1) \end{bmatrix} \in \R^{1 \times (N+1)}$,
\begin{subequations}\label{eq: components of aug truncate model}
\begin{align}
Y_a(t) = \begin{bmatrix} Y(t) \\ \zeta(t) \end{bmatrix} , \quad
Y_{\Phi,a}(\tau) = \begin{bmatrix} Y_\Phi(\tau) \\ \zeta_0(\tau) \end{bmatrix} , \\
A_a = \begin{bmatrix} A & 0 \\ C & 0 \end{bmatrix} , \quad
B_a = \begin{bmatrix} B \\ \alpha \end{bmatrix} , \quad
\Gamma(t) = \begin{bmatrix} 0 \\ - r(t) \end{bmatrix} .
\end{align}
\end{subequations}

\subsection{Proposed control strategy}

The proposed control strategy consists of a stabilizing state feedback of the truncated model (\ref{eq: aug truncate model}). Such a procedure is allowed by the following lemma. 

\begin{lem}\label{lem: Kalman condition}
$(A_a,B_a)$ satisfies the Kalman condition. 
\end{lem}

\textbf{Proof.}
We define the matrix $T = \begin{bmatrix} A & B \\ C & \alpha \end{bmatrix} \in \R^{(N+2) \times (N+2)}$. From (\ref{eq: components of aug truncate model}), the Hautus test shows that the pair $(A_a,B_a)$ satisfies the Kalman condition if an only if $T$ is invertible and the pair $(A,B)$ satisfies the Kalman condition. To show the former, let $Y_* = \begin{bmatrix} x_{*,n} & \ldots & x_{*,N} & u_* \end{bmatrix}^\top \in \operatorname{ker}(T)$. From (\ref{eq: components of truncate model}) and (\ref{eq: components of aug truncate model}) we deduce that $\lambda_n x_{*,n} + (a_n + \lambda_n b_n) u_* = 0$ for all $0 \leq n \leq N$ and $\sum_{n=0}^N x_{*,n} e_n(1) + \alpha u_* = 0$. Since $\lambda_{n} < 0$ for all $n \geq N+1$, we define $x_{*,n} = - \frac{a_n + \lambda_n b_n}{\lambda_n} u_*$ for all $n \geq N+1$. Hence we have $\lambda_n x_{*,n} + (a_n + \lambda_n b_n) u_* = 0$ for all $n \geq 0$. We also define $w_{*,n} = x_{*,n} + b_n u_*$ that gives $\lambda_n w_{*,n} + a_n u_* = 0$ for all $n \geq 0$. We infer that $(w_{*,n})_n$ and $(\lambda_n w_{*,n})_n$ are in $\ell^2(\N)$, hence we can define $w_* = \sum_{n \geq 0} w_{*,n} e_n \in D(\mathcal{A}_{c,0})$. Moreover, the latter equation shows that $\mathcal{A}_{c,0} w_* + \mathcal{A}_c L u_* = 0$. Using now the definition of $\alpha$ given by (\ref{eq: def alpha}), we deduce that $0 = \sum_{n=0}^N x_{*,n} e_n(1) + \alpha u_* = \sum_{n \geq 0} w_{*,n} e_n(1) = w_*(1)$. Therefore, introducing $x_* = w_* + L u_* \in D(\mathcal{A}_c)$, we obtain that $\mathcal{A}_c x_* = 0$ and $x_*(1) = w_*(1)+[Lu_*](1) = 0$. This shows that $a x_*'' + (b+c) x_* = 0$ with $x_*(1)=0$ and $x_*'(1) = -\cot(\theta) x_*(1) = 0$. So, by Cauchy uniqueness, $x_* = 0$. Since $w_*(0) = 0$, we get $0 = x_*(0) = [Lu_*](0) = u_*$. We infer $w_* = 0$ hence $w_{*,n} = 0$ for all $n \geq 0$. This implies that $x_{*,n} = w_{*,n} - b_n u_* = 0$ for all $n \geq 0$. We deduce that $Y_* = 0$, which shows that $T$ in invertible.

We now show that $(A,B)$ satisfies the Kalman condition. In view of (\ref{eq: components of truncate model}), since $A$ is diagonal with simple eigenvalues, it is sufficient to show that $a_n + \lambda_n b_n \neq 0$ for all $n \geq 0$. From (\ref{eq: def an bn}), using two integration by parts and the identity $\mathcal{A}_{c,0} e_n = \lambda_n e_n$, we obtain that $a_n + \lambda_n b_n = a e'_n(0)$. Since $e_n \neq 0$ with $e_n(0) = 0$, we obtain by Cauchy uniqueness that $e_n'(0) \neq 0$ hence $a_n + \lambda_n b_n \neq 0$ for all $n \geq 0$. Thus $(A,B)$ satisfies the Kalman condition, which completes the proof. 
\qed

Thus there exists $K \in \R^{1 \times (N+2)}$ such that $A_K = A_a + B_a K$ is Hurwitz with simple eigenvalues. We set for $t \geq 0$
\begin{equation}\label{eq: control input u}
u(t) = K Y_a(t) + p(t) 
\end{equation}
 where $p$ is a boundary disturbance. The control (\ref{eq: control input u}) takes the form of a PI controller because composed of 1) a proportional feedback of the state, via $Y(t)$, and 2) the integral component $\zeta(t)$ given by (\ref{eq: zeta dynamics}). We now need to select the integer $N \geq 0$ such that the closed-loop system composed of (\ref{eq: heat equation - ACS}), (\ref{eq: zeta dynamics}), and (\ref{eq: control input u}), is exponentially input-to-state stable with respect to the boundary perturbation $p$ and achieves the setpoint reference tracking of the system output $z(t) = y(t,1)$.

\subsection{Well-posedness of the closed-loop system dynamics}

The study of the well-posedness of the closed-loop system, which requires the introducing of the augmented state $X_\zeta = (X,\zeta)$ belonging to $\mathcal{H}_\zeta = L^2(0,1) \times \R$ endowed with the inner product $\left< (f,\zeta_f) , (g,\zeta_g) \right>_\zeta = \int_0^1 fg \,\mathrm{d}x + \zeta_f \zeta_g$, easily leads to the following result. 

\begin{lem}
Let $0 < h_m < h_M$, $h \in \mathcal{C}^0(\R_+)$ with $h_m \leq h(t) \leq h_M$, $\Phi \in \mathcal{C}^0([-h_M,0];\mathcal{H})$, $\zeta_0 \in \mathcal{C}^0([-h_M,0])$, $p \in \mathcal{C}^1(\R_+)$, and $r \in \mathcal{C}^0(\R_+)$. Then there exists a unique mild solution $X_\zeta = (X,\zeta) \in \mathcal{C}^0(\R_+;\mathcal{H}_\zeta)$ of (\ref{eq: heat equation - ACS}) and (\ref{eq: zeta dynamics}) with control input (\ref{eq: control input u}). Moreover we have $u,\zeta \in \mathcal{C}^1(\R_+)$.
\end{lem}

To assess the setpoint regulation, we need to resort to the concept of classical solutions. The existence and uniqueness of such solutions is guaranteed by the following corollary whose proof is an immediate consequence of classical results, see, e.g., \cite[Thm.~3.1.3]{curtain2012introduction}.

\begin{cor} 
Let $0 < h_m < h_M$, $h \in \mathcal{C}^1(\R_+)$ with $h_m \leq h(t) \leq h_M$ and such that $t \mapsto t-h(t)$ crosses $0$ a finite number of times, $\Phi \in \mathcal{C}^1([-h_M,0];\mathcal{H})$, $\zeta_0 \in \mathcal{C}^1([-h_M,0])$, $p \in \mathcal{C}^2(\R_+)$, and $r \in \mathcal{C}^1(\R_+)$. Assume that $\Phi(0) \in D(\mathcal{A})$ so that the compatibility condition
\begin{equation}\label{eq: classical solutions - compatibility condition}
\mathcal{B}\Phi(0) = K Y_{\Phi,a}(0) + p(0)
\end{equation}
holds. Then there exists a unique classical solution $X_\zeta = (X,\zeta) \in \mathcal{C}^0(\R_+;D(\mathcal{A})\times\R) \cap \mathcal{C}^1(\R_+;\mathcal{H}_\zeta)$ of (\ref{eq: heat equation - ACS}) and (\ref{eq: zeta dynamics}) with control input (\ref{eq: control input u}). Moreover we have $u \in \mathcal{C}_\mathrm{pw}^2(\R_+)$.
\end{cor}

\begin{rem}
From (\ref{eq: components of aug truncate model}) and because $A_K = A_a + B_a K$ is Hurwitz hence invertible, the last coefficient of $K$, that corresponds to the integral state $\zeta$, is necessarily non zero. Hence, for any given initial condition
$\Phi \in \mathcal{C}^1([-h_M,0];\mathcal{H})$ with $\Phi(0) \in D(\mathcal{A})$ and any boundary perturbation $p \in \mathcal{C}^2(\R_+)$, one can always select the initial condition $\zeta_0 \in \mathcal{C}^1([-h_M,0])$ of the integral component such that the compatibility condition (\ref{eq: classical solutions - compatibility condition}) holds.
\end{rem}

\section{Equilibrium conditions and associated dynamics of deviations}\label{sec: eq cond and dynamics of deviations}

\subsection{Equilibrium conditions}

Let $r_e,p_e \in\R$ be ``nominal'' values of the reference signal $r(t)$ and the boundary perturbation $p(t)$, respectively. Our first objective is to derive the equilibrium condition of the closed-loop system when setting $r(t) = r_e$ and $p(t)=p_e$. To do so, we denote by the subscript ``e'' the equilibrium condition associated with the different system signals. We define $Y_{a,e} = \begin{bmatrix} Y_e^\top & \zeta_e \end{bmatrix}^\top$, $Y_e = \begin{bmatrix} x_{0,e} & \ldots & x_{N,e} \end{bmatrix}^\top$, and $\Gamma_e = \begin{bmatrix} 0 & - r_e \end{bmatrix}^\top$. From (\ref{eq: aug truncate model}) and (\ref{eq: control input u}) we set $Y_{a,e} = - A_K^{-1} ( B_a p_e + \Gamma_e )$ and $u_e = K Y_{a,e} + p_e$ which give $0 = A_a Y_{a,e} + B_a u_e + \Gamma_e$. From (\ref{eq: components of truncate model}) and (\ref{eq: components of aug truncate model}), this implies that $0 = \lambda_n x_{n,e} + (a_n + \lambda_n b_n) u_e$ for all $0 \leq n \leq N$ and $0 = \sum_{n=0}^{N} x_{n,e} e_n(1) + \alpha u_e - r_e$. Regarding the residual dynamics given by (\ref{eq: ODE coefficient of projection}) for $n \geq N+1$, we define $x_{n,e} = -\frac{a_n + \lambda_n b_n}{\lambda_n} u_e$. This yields $0 = \lambda_n x_{n,e} + (a_n + \lambda_n b_n) u_e$ for all $n \geq 0$. We note that $(x_{n,e})_{n \geq 0} \in \ell^2(\mathbb{N})$ hence we can define $X_e = \sum_{n \geq 0} x_{n,e} e_n \in \mathcal{H}$. Moreover, introducing for $n \geq 0$ the quantities $w_{n,e} = x_{n,e} + b_n u_e$, we have for $n \geq N+1$ that $w_{n,e} = -\frac{a_n}{\lambda_n} u_e$, showing that $(w_{n,e})_{n \geq 0} \in \ell^2(\mathbb{N})$ and $(\lambda_n w_{n,e})_{n \geq 0} \in \ell^2(\mathbb{N})$. This allows the introduction of $W_e = \sum_{n \geq 0} w_{n,e} e_n \in D(\mathcal{A}_0) = D(\mathcal{A}_{c,0})$. Moreover, from the definition of $b_n$ given by (\ref{eq: def an bn}), we have $X_e = W_e + L u_e \in D(\mathcal{A}_c)$ hence $\mathcal{B}X_e = u_e$. Furthermore, since $\lambda_n w_{n,e} + a_n u_e = 0$ for all $n \geq 0$, we have from the definition of $a_n$ given by (\ref{eq: def an bn}) that $\mathcal{A}_{c,0} W_e + \mathcal{A}_{c}Lu_e = 0$ hence $\mathcal{A}_{c}X_e = 0$. Using now Lemma~\ref{lem: seies expansion of f(1)}, (\ref{eq: def alpha}), and the above relations between $x_{n,e}$ and $w_{n,e}$, we obtain from $0 = \sum_{n=0}^{N} x_{n,e} e_n(1) + \alpha u_e - r_e$ that $W_e(1) = r_e$. Since $X_e \in D(\mathcal{A}_c) \subset H^1(0,1)$, we infer that $X_e(1) = W_e(1) + [Lu_e](1) = r_e$, which provides the desired reference tracking.

\subsection{Dynamics of deviations}

Let $r_e,p_e \in \R$ be arbitrary and consider the different equilibrium quantities defined above. We can introduce the dynamics of deviations of the system trajectory with respect to the considered equilibrium condition. These deviations are denoted by the symbol ``$\Delta$''. For instance, $\Delta X(t)$ stands for $X(t)-X_e$. We obtain the following dynamics of deviation: $\frac{\mathrm{d}(\Delta X)}{\mathrm{d}t}(t) = \mathcal{A} \Delta X(t) + c \Delta  X(t-h(t))$, $\mathcal{B} \Delta X(t) = \Delta u(t)$, $\Delta \dot{\zeta}(t) = \sum_{n = 0}^N \Delta x_n(t) e_n(1) + c \{ \Delta \zeta(t-h(t)) - \Delta \zeta(t) \} + \alpha \Delta u(t) - \Delta r(t)$, $\Delta x_n(t) = \langle \Delta X(t) , e_n \rangle$ and $\Delta w_n(t) = \langle \Delta W(t) , e_n \rangle = \Delta x_n(t) + b_n \Delta u(t)$. This yields the following representation for the closed-loop system dynamics:
\begin{subequations}\label{eq: dynamics of deviations}
\begin{align}
\Delta \dot{Y}_a(t) & = A_K \Delta Y_a(t) + c \{ \Delta Y_a(t-h(t)) - \Delta Y_a(t) \}  \nonumber \\
& \phantom{=}\; + B_a \Delta p(t) + \Delta \Gamma(t) , \label{eq: dynamics of deviations - truncated model} \\
\Delta \dot{x}_n(t) & = \lambda_n \Delta x_n(t) + c \{ \Delta x_n(t-h(t)) - \Delta x_n(t) \} \nonumber \\
& \phantom{=}\; + (a_n + \lambda_n b_n) \Delta u(t) , \quad n \geq N+1 , \label{eq: dynamics of deviations - residual dynamics} \\
\Delta u(t) & = K \Delta Y_a(t) + \Delta p(t) \label{eq: dynamics of deviations - control input}\\
\Delta Y_a(\tau) & = \Delta Y_{\Phi,a}(\tau) , \quad \tau\in[-h_M,0] \label{eq: dynamics of deviations - IC trunated model} \\
\Delta x_n(\tau) & = \langle \Delta \Phi(\tau) , e_n \rangle , \quad \tau\in[-h_M,0] , \quad n \geq 0 \label{eq: dynamics of deviations - IC residual dynamics}
\end{align}
\end{subequations}

\section{Stability analysis}\label{sec: stability analysis}

The main result of this section is stated as follows.

\begin{thm}\label{thm: stability}
Let $0 < h_m < h_M$ be arbitrarily given. Let $N \geq 0$ be such that $\lambda_{N + 1} < - 2 \sqrt{5} \vert c \vert$ and consider the matrices $A_a$ and $B_a$ defined by (\ref{eq: components of aug truncate model}). Let $K \in \mathbb{R}^{1 \times (N+2)}$ be such that $A_K = A_a + B_a K$ is Hurwitz with simple eigenvalues $\mu_1 , \ldots , \mu_{N+2} \in \mathbb{C}$ satisfying $\operatorname{Re}\mu_n < - 3 \vert c \vert$ for all $1 \leq n \leq N+2$. Then, there exist constants $\kappa , \overline{C}_0 , \overline{C}_1 > 0$ such that, for all $h \in \mathcal{C}^0(\R_+)$ with $h_m \leq h(t) \leq h_M$, $\Phi \in \mathcal{C}^0([-h_M,0];\mathcal{H})$, $\zeta_0 \in \mathcal{C}^0([-h_M,0])$, $p \in \mathcal{C}^1(\R_+)$, and $r \in \mathcal{C}^0(\R_+)$, the mild solution $X_\zeta = (X,\zeta) \in\mathcal{C}^0(\mathbb{R}_+;\mathcal{H}_\zeta)$ of (\ref{eq: heat equation - ACS}) and (\ref{eq: zeta dynamics}) with control input (\ref{eq: control input u}) satisfies, for all $t \geq 0$,
\begin{align}
& \Vert \Delta X(t) \Vert + \vert \Delta \zeta(t) \vert + \vert \Delta u(t) \vert \nonumber \\
& \quad \leq \overline{C}_0 e^{- \kappa t} \sup\limits_{\tau \in [-h_M,0]} ( \Vert \Delta \Phi(\tau) \Vert + \vert \Delta \zeta_0(\tau) \vert ) \label{eq: main theorem - estimate trajectory}\\
& \quad \phantom{\leq}\; + \overline{C}_1 \sup\limits_{\tau \in [0,t]} e^{- \kappa (t-\tau)} \left( \vert \Delta p(\tau) \vert + \vert \Delta r(\tau) \vert \right) \nonumber .
\end{align}
\end{thm}

\begin{cor}
In the context of Theorem~\ref{thm: stability}, assume that $r(t) \rightarrow r_e$ and $p(t) \rightarrow p_e$ as $t \rightarrow + \infty$. Then $X(t) \rightarrow X_e$ and $\zeta(t) \rightarrow \zeta_e$ as $t \rightarrow + \infty$ with exponential vanishing of the contribution of the initial conditions. 
\end{cor}

\begin{rem}
From Theorem~1, one needs to start by selecting the integer $N \geq 0$ such that $\lambda_{N + 1} < - 2 \sqrt{5} \vert c \vert$. This is always possible because $\lambda_n \sim -an^2\pi^2$ as $n \rightarrow +\infty$ with $a>0$. Then, because of Lemma~\ref{lem: Kalman condition}, the feedback gain $K \in \mathbb{R}^{1 \times (N+2)}$ can always be computed such that $A_K = A_a + B_a K$ is Hurwitz with arbitrary eigenvalue assignment. This allows the application of Theorem~\ref{thm: stability}.
\end{rem}

\subsection{Truncated model}\label{subsec: stability truncated model}

The design of the feedback gain $K$ and the resulting stability properties of the truncated model (\ref{eq: dynamics of deviations - truncated model}) rely on the following lemma whose proof is identical to~\cite[Lem.~8]{lhachemi2020boundary}.

\begin{lem}\label{lem: prel lemma truncated model}
Let $N \geq 1$, $0 < h_m < h_M$, $\mathbf{A} \in \R^{N \times N}$, and $c \in \R$. Assume that $\mathbf{A}$ is Hurwitz with simple eigenvalues $\mu_1,\ldots,\mu_N \in \C$ such that $\operatorname{Re}\mu_n < -3 \vert c \vert$ for all $1 \leq n \leq N$. Then there exist $\sigma,C_0,C_1 > 0$ such that, for any $x_0 \in \mathcal{C}^0([-h_M,0];\R^N)$, any $h \in \mathcal{C}^0(\R_+)$ with $h_m \leq h(t) \leq h_M$, and any $q \in \mathcal{C}^0(\R_+;\R^N)$, the trajectory of
\begin{align*}
\dot{x}(t) & = \mathbf{A} x(t) + c \left\{ x(t-h(t)) - x(t) \right\} + q(t) , \\
x(\tau) & = x_0(\tau) , \quad \tau \in [-h_M,0]
\end{align*}
satisfies, for all $t \geq 0$,
\begin{align}
\Vert x(t) \Vert 
& \leq C_0 e^{-\sigma t} \sup\limits_{\tau\in[-h_M,0]} \Vert x_0(\tau) \Vert \label{eq: prel lemma truncated model - ISS} \\
& \phantom{\leq}\; + C_1 \sup\limits_{\tau\in[0,t]} e^{-\sigma (t-\tau)} \Vert q(\tau) \Vert . \nonumber
\end{align}
\end{lem}
From the assumptions of Thm.~\ref{thm: stability}, Lemma~\ref{lem: prel lemma truncated model} applies to the truncated model (\ref{eq: dynamics of deviations - truncated model}) with initial condition (\ref{eq: dynamics of deviations - IC trunated model}).

\subsection{Residual infinite-dimensional dynamics}\label{subsec: stability residual dynamics}

We now need to investigate the selection of the integer $N \geq 0$ such that the residual dynamics composed of (\ref{eq: dynamics of deviations - residual dynamics}) and (\ref{eq: dynamics of deviations - IC residual dynamics}) is exponentially stable. 

\begin{lem}\label{lem: infinite-dim part negelected in the design}
Let $0 < h_m < h_M$ and $\sigma,C_2,C_3>0$ be arbitrarily given. Let $N \geq 0$ be such that $\lambda_{N + 1} < - 2 \sqrt{5} \vert c \vert$. Then, there exist constants $\kappa\in(0,\sigma)$ and $C_4 , C_5 > 0$ such that, for all $h \in \mathcal{C}^0(\R_+)$ with $h_m \leq h(t) \leq h_M$, $\Phi \in \mathcal{C}^0([-h_M,0];\mathcal{H})$, $\zeta_0 \in \mathcal{C}^0([-h_M,0])$, $p \in \mathcal{C}^1(\R_+)$, $r \in \mathcal{C}^0(\R_+)$, and $u \in \mathcal{C}^1(\mathbb{R}_+)$ such that
\begin{align}
\vert \Delta u(t) \vert 
& \leq C_2 e^{-\sigma t} \sup\limits_{\tau\in[-h_M,0]} \left( \Vert \Delta \Phi(\tau) \Vert + \vert \Delta \zeta_0(\tau) \vert \right) \nonumber \\
& \phantom{\leq}\; + C_3 \sup\limits_{\tau\in[0,t]} e^{-\sigma (t-\tau)} \left( \vert \Delta p(\tau) \vert + \vert \Delta r(\tau) \vert \right) \label{eq: assumed estimate on u}
\end{align}
for all $t \geq 0$, the mild solution $X_\zeta = (X,\zeta) \in\mathcal{C}^0(\mathbb{R}_+;\mathcal{H}_\zeta)$ of (\ref{eq: heat equation - ACS}) and (\ref{eq: zeta dynamics}) satisfies for all $t \geq 0$
\begin{align}
& \sum\limits_{n \geq N + 1} \vert \Delta x_n(t) \vert^2 \nonumber \\
& \leq C_4 e^{- 2 \kappa t} \sup\limits_{\tau \in [-h_M,0]} \left( \Vert \Delta \Phi(\tau) \Vert + \vert \Delta \zeta_0(\tau) \vert \right)^2 \nonumber \\
& \phantom{\leq}\; + C_5 \sup\limits_{\tau \in [0,t]} e^{- 2 \kappa (t-\tau)} \left( \vert \Delta p(\tau) \vert + \vert \Delta r(\tau) \vert \right)^2 . \label{eq: exp stab infinite-dim part}
\end{align}
\end{lem}

\begin{rem}
The design constraint $\lambda_{N + 1} < - 2 \sqrt{5} \vert c \vert$ is the same as in~\cite[Lem.~10]{lhachemi2020boundary}. However, the proof reported therein does not apply in the presence of the boundary perturbation $p$. Indeed, following the lines of~\cite[Lem.~10]{lhachemi2020boundary}, one gets an estimate similar to (\ref{eq: exp stab infinite-dim part}) but with the occurrence of the extra term $\vert \Delta\dot{p}(\tau) \vert$ in the term evaluating the contribution of $\Delta p$ and $\Delta r$. We refine here the stability analysis in order to obtain the claimed estimate (\ref{eq: exp stab infinite-dim part}) involving only $\Delta p$, and not $\Delta\dot{p}$. 
\end{rem}

\textbf{Proof.}
Let $N \geq 0$ be such that $\lambda_{N + 1} < - 2 \sqrt{5} \vert c \vert$. We define $\eta = - \lambda_{N+1}/2 > \sqrt{5} \vert c \vert \geq 0$, which is such that $\lambda_n \leq \lambda_{N+1} = -2\eta <0$ for all $n \geq N+1$. Note that, in this proof, we always consider integers $n \geq N+1$. Let $\kappa \in ( 0 , \min(\eta,\sigma) )$ be arbitrarily given and to be specified later. We introduce, for $t \geq 0$, $\Delta v_n(t) = \Delta x_n(t) - \Delta x_n(t-h(t))$, yielding
\begin{equation}\label{eq: dynamics of deviation xn - funcion vn}
\Delta \dot{x}_n(t) = \lambda_n \Delta x_n(t) - c \Delta v_n(t) + (a_n + \lambda_n b_n) \Delta u(t)
\end{equation}
for all $t \geq 0$. We also consider the series
\begin{align*}
S_x(t) & = \sum_{n \geq N+1} \vert \Delta x_n(t) \vert^2 , & t \geq -h_M ; \\
S_v(t) & = \sum_{n \geq N+1} \vert \Delta v_n(t) \vert^2 , & t \geq 0 
\end{align*}
which are finite because $S_x(t) \leq \Vert \Delta X(t) \Vert^2$ and $S_v(t) \leq 2 S_x(t) + 2 S_x(t-h(t))$. Finally, we introduce for any $t_1 < t_2$ and any real-valued and continuous function $\psi$ the notation $\mathcal{I}(\psi,t_1,t_2) = \int_{t_1}^{t_2} e^{- 2 \eta (t_2-\tau)} \vert \psi(\tau) \vert \,\mathrm{d}\tau$. We have $\mathcal{I}(\psi,t_1,t_2) \leq \frac{1 - e^{-2(\eta-\kappa)(t_2-t_1)}}{2 (\eta-\kappa)} \sup\limits_{\tau\in[t_1,t_2]} e^{-2\kappa(t_2-\tau)} \vert \psi(\tau) \vert$ and $\mathcal{I}(\psi,t_1,t_2)^2 \leq \frac{1 - e^{-2\eta(t_2-t_1)}}{2 \eta} \mathcal{I}(\psi^2,t_1,t_2)$. By integrating (\ref{eq: dynamics of deviation xn - funcion vn}), we obtain for $t \geq h_M$
\begin{align*}
& \Delta v_n(t) = \left\{ e^{\lambda_n h(t)} - 1 \right\} \Delta x_n(t-h(t)) \\
& + \int_{t-h(t)}^{t} e^{\lambda_n(t-\tau)} \{ - c \Delta v_n(\tau) + (a_n+\lambda_n b_n) \Delta u(\tau) \} \,\mathrm{d}\tau
\end{align*}
hence, using $\lambda_n \leq - 2 \eta$,
\begin{align*}
& \vert \Delta v_n(t) \vert 
\leq \vert \Delta x_n(t-h(t)) \vert 
+ \vert c \vert \mathcal{I}(\Delta v_n,t-h(t),t) \\
& \qquad + \vert a_n \vert \mathcal{I}(\Delta u,t-h(t),t) \\
& \qquad + \vert b_n \vert \left\vert \lambda_n \int_{t-h(t)}^{t} e^{\lambda_n(t-\tau)} \Delta u(\tau) \,\mathrm{d}\tau \right\vert .
\end{align*}
Since $\kappa < \eta$ we have $\left\vert \lambda_n \int_{t-h(t)}^{t} e^{\lambda_n(t-\tau)} \Delta u(\tau) \,\mathrm{d}\tau \right\vert \leq \dfrac{2\eta}{2 \eta - \kappa} \sup_{\tau\in[t-h(t),t]} e^{-\kappa(t-\tau)} \vert \Delta u(\tau) \vert$ because $\lambda_n \leq - 2\eta < -\eta < - \kappa < 0$. Combining the two latter estimates and using Young's inequality we obtain
\begin{align*}
& \vert \Delta v_n(t) \vert^2
\leq 4 \vert \Delta x_n(t-h(t)) \vert^2 
+ \gamma_1 \vert c \vert^2 \mathcal{I}(\Delta v_n^2,t-h(t),t) \\ 
& \qquad + \gamma_1 \vert a_n \vert^2 \mathcal{I}(\Delta u^2,t-h(t),t) \\
& \qquad + \dfrac{16\eta^2}{(2\eta-\kappa)^2} \vert b_n \vert^2 \sup_{\tau\in[t-h(t),t]} e^{-2\kappa(t-\tau)} \vert \Delta u(\tau) \vert^2 
\end{align*}
for all $t \geq h_M$ where $\gamma_1 = \frac{2}{\eta} (1-e^{-2\eta h_M})$. Summing for $n \geq N+1$, we obtain for $t \geq h_M$
\begin{align*}
S_v(t)
& \leq 4 S_x(t-h(t))
+ \gamma_2(\kappa) \vert c \vert^2 \sup_{\tau\in[t-h(t),t]} e^{-2\kappa(t-\tau)} S_v(\tau) \\ 
& \phantom{\leq}\; + \gamma_3(\kappa) \sup_{\tau\in[t-h(t),t]} e^{-2\kappa(t-\tau)} \vert \Delta u(\tau) \vert^2
\end{align*}
where $a = \mathcal{A}_{c} L \mathbbm{1}$, $b = - L \mathbbm{1}$, $\gamma_2(\kappa) = \frac{1}{\eta(\eta-\kappa)} (1-e^{-2\eta h_M}) (1-e^{-2(\eta-\kappa) h_M})$ and $\gamma_3(\kappa) = \gamma_2(\kappa) \Vert a \Vert^2 + \frac{16\eta^2}{(2\eta-\kappa)^2} \Vert b \Vert^2$. This implies that, for all $t \geq h_M$,
\begin{align}
& \sup_{\tau\in[h_M,t]} e^{2\kappa\tau} S_v(\tau) 
\leq 4 e^{2\kappa h_M} \sup_{\tau\in[0,t-h_m]} e^{2\kappa\tau} S_x(\tau) \label{eq: stability residual dynamics - prel 1} \\
& + \gamma_2(\kappa) \vert c \vert^2 \sup_{\tau\in[0,t]} e^{2\kappa\tau} S_v(\tau) 
+ \gamma_3(\kappa) \sup_{\tau\in[0,t]} e^{2\kappa\tau} \vert \Delta u(\tau) \vert^2 \nonumber . 
\end{align}

Integrating now (\ref{eq: dynamics of deviation xn - funcion vn}) on $[0,t]$ for $t \geq 0$, using again $\lambda_n \leq - 2\eta$, and proceeding as in the previous paragraph, we infer that, for all $t \geq 0$,
\begin{align}
S_x(t) 
& \leq 4 e^{-2\kappa t} S_x(0) + \gamma_4(\kappa) \vert c \vert^2 \sup_{\tau \in [0,t]} e^{-2\kappa(t-\tau)} S_v(\tau) \nonumber \\
& \phantom{\leq}\;  + \gamma_5(\kappa) \sup_{\tau\in[0,t]} e^{-2\kappa(t-\tau)} \vert \Delta u(\tau) \vert^2 \label{eq: stability residual dynamics - prel 2} 
\end{align}
where $\gamma_4(\kappa) = \frac{1}{\eta(\eta-\kappa)}$ and $\gamma_5(\kappa) = \gamma_4(\kappa) \Vert a \Vert^2 + \frac{16\eta^2}{(2\eta-\kappa)^2} \Vert b \Vert^2$. Combining (\ref{eq: stability residual dynamics - prel 1}-\ref{eq: stability residual dynamics - prel 2}) and noting that $S_x(0) \leq \Vert \Delta \Phi(0) \Vert^2$, we obtain for $t \geq h_M$
\begin{align*}
& \sup_{\tau\in[h_M,t]} e^{2\kappa\tau} S_v(\tau) 
\leq 16 e^{2\kappa h_M} \Vert \Delta \Phi(0) \Vert^2 \\
& \qquad + \xi(\kappa) \sup_{\tau\in[0,t]} e^{2\kappa\tau} S_v(\tau) 
+ \gamma_6(\kappa) \sup_{\tau\in[0,t]} e^{2\kappa\tau} \vert \Delta u(\tau) \vert^2
\end{align*}
with $\gamma_6(\kappa) = \gamma_3(\kappa) + 4 e^{2\kappa h_M} \gamma_5(\kappa)$ and
\begin{align*}
& \xi(\kappa)
= \gamma_2(\kappa) \vert c \vert^2 + 4 e^{2\kappa h_M} \gamma_4(\kappa) \vert c \vert^2 \\
& = \dfrac{\vert c \vert^2}{\eta(\eta-\kappa)} \left\{ 4 e^{2\kappa h_M} + (1-e^{-2\eta h_M}) (1-e^{-2(\eta-\kappa) h_M}) \right\} .
\end{align*}
Recalling from the design constraint $\lambda_{N + 1} < - 2 \sqrt{5} \vert c \vert$ that $\eta > \sqrt{5} \vert c \vert$, we have $5 \vert c \vert^2 / \eta^2 < 1$. Hence, a continuity argument at $\kappa = 0$ shows the existence of $\kappa \in ( 0 , \min(\eta,\sigma) )$ such that $0 \leq \xi(\kappa) < 1$. We fix such a $\kappa \in ( 0 , \min(\eta,\sigma) )$ for the rest of the proof. Since all the considered supremums are finite, we deduce from the latter estimate that, for all $t \geq h_M$,
\begin{align}
& \sup_{\tau\in[h_M,t]} e^{2\kappa\tau} S_v(\tau) 
\leq \dfrac{16 e^{2\kappa h_M}}{1-\xi} \Vert \Delta \Phi(0) \Vert^2 \label{eq: stability residual dynamics - prel 3} \\
& + \dfrac{\xi}{1-\xi} \sup_{\tau\in[0,h_M]} e^{2\kappa\tau} S_v(\tau) + \dfrac{\gamma_6}{1-\xi} \sup_{\tau\in[0,t]} e^{2\kappa\tau} \vert \Delta u(\tau) \vert^2  \nonumber
\end{align}
where we dropped the dependency of $\gamma_6,\xi$ on the parameter $\kappa$ which is now fixed. To conclude the proof, we need to estimate the term $\sup_{\tau\in[0,t]} e^{2\kappa\tau} S_v(\tau)$ for $t \in [0,h_M]$. From the definition of $S_v$ we have, for any $t \in [0,h_M]$, $\sup_{\tau\in[0,t]} e^{2\kappa\tau} S_v(\tau) \leq 4 e^{2\kappa h_M} \sup_{\tau\in[-h_M,t]} S_x(\tau)$. From (\ref{eq: dynamics of deviations - residual dynamics}) and recalling that $n \geq N+1$ with $\lambda_n \leq - 2 \eta < - 2 \sqrt{5} \vert c \vert$, we have $\lambda_n - c \leq \lambda_n + \vert c \vert < - (2 \sqrt{5} -1) \vert c \vert \leq 0$ hence
\begin{align*}
\vert \Delta x_n(t) \vert
& \leq \vert \Delta x_n(0) \vert + \vert c \vert \sqrt{h_M} \sqrt{\int_0^t  \vert \Delta x_n(\tau-h(\tau)) \vert^2 \,\mathrm{d}\tau} \\
& \phantom{\leq}\; + ( \vert a_n \vert h_M + \vert b_n \vert e^{\vert c \vert h_M} ) \sup_{\tau\in[0,t]} \vert \Delta u(\tau) \vert 
\end{align*}
for all $t \in [0,h_M]$. Using Young's inequality and summing for $n \geq N+1$, we obtain
\begin{align*}
S_x(t)
& \leq 3 S_x(0) + 3 \vert c \vert^2 h_M^2 \sup_{\tau\in[-h_M,t-h_m]} S_x(\tau)  \\
& \phantom{\leq}\; + 6 ( \Vert a \Vert^2 h_M^2 + \Vert b \Vert^2 e^{2 \vert c \vert h_M} ) \sup_{\tau\in[0,t]} \vert \Delta u(\tau) \vert^2
\end{align*}
for all $t \in [0,h_M]$. This implies, for all $t \in [0,h_M]$,
\begin{align*}
& \sup_{\tau\in[0,t]} S_x(\tau)
\leq 3 (1+\vert c \vert^2 h_M^2) \sup_{\tau\in[-h_M,0]} \Vert \Delta\Phi(\tau) \Vert^2 \\
& \qquad + 3 \vert c \vert^2 h_M^2 \sup_{\tau\in[0,\max(t-h_m,0)]} S_x(\tau) \\
& \qquad + 6 ( \Vert a \Vert^2 h_M^2 + \Vert b \Vert^2 e^{2 \vert c \vert h_M} ) \sup_{\tau\in[0,t]} \vert \Delta u(\tau) \vert^2 .
\end{align*}
By a simple induction argument (since $h_m > 0$), we obtain the existence of a constant $\gamma_7 > 0$ such that, for all $t \in [0,h_M]$, $\sup_{\tau\in[0,t]} S_x(\tau) \leq \gamma_7 \sup_{\tau\in[-h_M,0]} \Vert \Delta\Phi(\tau) \Vert^2 + \gamma_7 \sup_{\tau\in[0,t]} \vert \Delta u(\tau) \vert^2$. We deduce (see beginning of this paragraph) the existence of a constant $\gamma_8 > 0$ such that, for all $t \in [0,h_M]$, $\sup_{\tau\in[0,t]} e^{2\kappa\tau} S_v(\tau) 
\leq \gamma_8 \sup_{\tau\in[-h_M,0]} \Vert \Delta\Phi(\tau) \Vert^2 + \gamma_8 \sup_{\tau\in[0,t]} e^{2\kappa\tau} \vert \Delta u(\tau) \vert^2$. Combining this latter estimate with (\ref{eq: stability residual dynamics - prel 3}), we infer the existence of a constant $\gamma_9 > 0$ such that, for all $t \geq 0$, 
\begin{align}
\sup_{\tau\in[0,t]} e^{2\kappa\tau} S_v(\tau) 
& \leq \gamma_9 \sup_{\tau\in[-h_M,0]} \Vert \Delta\Phi(\tau) \Vert^2 \label{eq: stability residual dynamics - prel 4} \\
& \phantom{\leq}\; + \gamma_9 \sup_{\tau\in[0,t]} e^{2\kappa\tau} \vert \Delta u(\tau) \vert^2 . \nonumber
\end{align}
Substituting this estimate into (\ref{eq: stability residual dynamics - prel 2}), we obtain the existence of a constant $\gamma_{10} > 0$ such that, for all $t \geq 0$, $S_x(t) \leq \gamma_{10} e^{-2\kappa t} \sup_{\tau\in[-h_M,0]} \Vert \Delta\Phi(\tau) \Vert^2 + \gamma_{10} \sup_{\tau\in[0,t]} e^{-2\kappa(t-\tau)} \vert \Delta u(\tau) \vert^2$.
The claimed estimate (\ref{eq: exp stab infinite-dim part}) now directly follows from the assumption that $u$ satisfies (\ref{eq: assumed estimate on u}) and the fact that $0 < \kappa < \sigma$.
\qed

\subsection{Completion of the proof of Theorem~\ref{thm: stability}}

By applying first the result of Subsection~\ref{subsec: stability truncated model} and then the result of Subsection~\ref{subsec: stability residual dynamics}, the claimed estimate (\ref{eq: main theorem - estimate trajectory}) follows from $\vert \Delta \zeta(t) \vert \leq \Vert \Delta Y_a(t) \Vert$, $\Vert \Delta X(t) \Vert \leq \Vert \Delta Y_a(t) \Vert + \sqrt{\sum_{n \geq N + 1} \vert \Delta x_n(t) \vert^2}$, and (\ref{eq: dynamics of deviations - control input}). This completes the proof of Theorem~\ref{thm: stability}.

\section{Setpoint regulation assessment}\label{sec: setpoint regulation assessment}

We now address the setpoint regulation of the closed-loop system for classical solutions.

\begin{thm}\label{thm: regulation}
Under the assumptions of Theorem~\ref{thm: stability}, and for the same constant $\kappa > 0$, there exist constants $\overline{C}_2 , \overline{C}_3 > 0$ such that, for all $h \in \mathcal{C}^1(\R_+)$ with $h_m \leq h(t) \leq h_M$ and so that $t \mapsto t-h(t)$ crosses $0$ a finite number of times, $\Phi \in \mathcal{C}^1([-h_M,0];\mathcal{H})$ with $\Phi(0) \in D(\mathcal{A})$, $\zeta_0 \in \mathcal{C}^1([-h_M,0])$, $p \in \mathcal{C}^2(\R_+)$, and $r \in \mathcal{C}^1(\R_+)$, all such that the compatibility condition (\ref{eq: classical solutions - compatibility condition}) holds, we have, for all $t \geq 0$,
\begin{align}
& \vert [X(t)](1) - r(t) \vert \leq \label{eq: main theorem - estimate regulation} \\
& \overline{C}_2 e^{- \kappa t} \left\{ \sup\limits_{\tau \in [-h_M,0]} ( \Vert \Delta \Phi(\tau) \Vert + \vert \Delta \zeta_0(\tau) \vert ) + \Vert \mathcal{A}_c \Delta\Phi(0) \Vert \right\} \nonumber \\
& + \overline{C}_3 \sup\limits_{\tau \in [0,t]} e^{- \kappa (t-\tau)} \left( \vert \Delta p(\tau) \vert + \vert \Delta \dot{p}(\tau) \vert + \vert \Delta r(\tau) \vert \right) \nonumber
\end{align}
\end{thm}

\begin{cor}\label{cor: regulation}
In the context of Theorem~\ref{thm: regulation}, assume that $r(t) \rightarrow r_e$, $p(t) \rightarrow p_e$, and $\dot{p}(t) \rightarrow 0$ as $t \rightarrow + \infty$. Then $[X(t)](1) \rightarrow r_e$ as $t \rightarrow + \infty$ with exponential vanishing of the contribution of the initial conditions. 
\end{cor}

\textbf{Proof of Theorem~\ref{thm: regulation}.}
Recalling that, for classical solutions, $W(t) = X(t) - Lu(t) \in D(\mathcal{A}_{c,0})$, and since $W_e = X_e - L u_e \in D(\mathcal{A}_{c,0})$ with $X_e(1) = W_e(1) = r_e$, we have $\vert [X(t)](1) - r(t) \vert \leq \vert [W(t)](1) - r_e \vert + \vert \Delta r(t) \vert \leq \vert [\Delta W(t)](1) \vert + \vert \Delta r(t) \vert$. To obtain (\ref{eq: main theorem - estimate regulation}), we only need to investigate the term $\vert [\Delta W(t)](1) \vert$. To do so, since $\lambda_n \sim -an^2\pi^2$ as $n \rightarrow +\infty$, let $\delta > 0$ and an integer $M \geq N$ be such that $\lambda_n \leq - (2\kappa+\delta) < 0$ and $\vert \lambda_n \vert \leq \vert \lambda_n \vert^2$ for all $n \geq M+1$. Then we have $\vert [\Delta W(t)](1) \vert = \left\vert \sum_{n \geq 0} \Delta w_n(t) e_n(1) \right\vert \leq \sqrt{\sum_{n = 0}^{M} \vert e_n(1) \vert^2} \Vert \Delta W(t) \Vert + \sqrt{\sum_{n  \geq M+1} \frac{ \vert e_n(1) \vert^2}{\vert \lambda_n \vert}} \sqrt{\sum_{n  \geq M+1} \vert \lambda_n \vert \vert \Delta w_n(t) \vert^2}$ where it can be seen from $\lambda_n \sim -an^2\pi^2$ and $e_n(1) = O(1)$ as $n \rightarrow +\infty$ that $\sum_{n  \geq M+1} \frac{ \vert e_n(1) \vert^2}{\vert \lambda_n \vert} < \infty$. From Theorem~\ref{thm: stability} and since $\Vert \Delta W(t) \Vert \leq \Vert \Delta X(t) \Vert + \Vert L \Vert \vert \Delta u(t) \vert$, we only need to study the term $\sum_{n  \geq M+1} \vert \lambda_n \vert \vert \Delta w_n(t) \vert^2$ to conclude that (\ref{eq: main theorem - estimate regulation}) holds. In the sequel we always consider integers $n \geq M+1$. From (\ref{eq: dynamics of deviations - residual dynamics}) and recalling that $\Delta w_n(t) = \Delta x_n(t) + b_n \Delta u(t)$, we have for $t \geq 0$ that $\Delta \dot{w}_n(t) = \lambda_n \Delta w_n(t) - c \Delta v_n(t) + a_n \Delta u(t) + b_n \Delta \dot{u}(t)$  with $\Delta v_n(t) = \Delta x_n(t) - \Delta x_n(t-h(t))$. Then we obtain after integration on $[0,t]$ that $\sqrt{\vert \lambda_n \vert} \vert \Delta w_n(t) \vert \leq e^{\lambda_n t} \sqrt{\vert \lambda_n \vert} \vert \Delta w_n(0) \vert + \vert c \vert \mathcal{J}_{1,n}(t) + \vert a_n \vert \mathcal{J}_{2,n}(t) + \vert b_n \vert \mathcal{J}_{3,n}(t)$ for all $t \geq 0$ and $n \geq M+1$ with $\mathcal{J}_{1,n}(t) = \sqrt{\vert \lambda_n \vert} \int_0^t e^{\lambda_n(t-\tau)} \vert \Delta v_n(\tau) \vert \,\mathrm{d}\tau$, $\mathcal{J}_{2,n}(t) = \sqrt{\vert \lambda_n \vert} \int_0^t e^{\lambda_n(t-\tau)} \vert \Delta u(\tau) \vert \,\mathrm{d}\tau$, and $\mathcal{J}_{3,n}(t) = \sqrt{\vert \lambda_n \vert} \int_0^t e^{\lambda_n(t-\tau)} \vert \Delta \dot{u}(\tau) \vert \,\mathrm{d}\tau$. Using $\lambda_{n} \leq -(2\kappa+\delta)$ and $\vert \lambda_n \vert \leq \vert \lambda_n \vert^2$ for all $n \geq M+1$, we obtain that $\mathcal{J}_{1,n}(t) \leq \sqrt{\int_0^t e^{-(2\kappa+\delta)(t-\tau)} \vert \Delta v_n(\tau) \vert^2 \,\mathrm{d}\tau}$, $\mathcal{J}_{2,n}(t) \leq 2 \sup_{\tau\in[0,t]} e^{-\kappa(t-\tau)} \vert \Delta u(\tau) \vert$, and $\mathcal{J}_{3,n}(t) \leq 2 \sup_{\tau\in[0,t]} e^{-\kappa(t-\tau)} \vert \Delta \dot{u}(\tau) \vert$. Combining the four latter inequalities, using next Young's inequality, and finally summing for $n \geq M+1$, we obtain that
\begin{align*}
& \sum_{n \geq M+1} \vert \lambda_n \vert \vert \Delta w_n(t) \vert^2
\leq 4 e^{- 2 \kappa t} \sum_{n \geq M+1} \vert \lambda_n \vert \vert \Delta w_n(0) \vert^2 \\
& \phantom{\leq}\; + 4 \vert c \vert^2 \int_0^t e^{-(2\kappa+\delta)(t-\tau)} \sum_{n \geq M+1} \vert \Delta v_n(\tau) \vert^2 \,\mathrm{d}\tau \\
& \phantom{\leq}\; + 16 \Vert a \Vert^2 \sup_{\tau\in[0,t]} e^{-2\kappa(t-\tau)} \vert \Delta u(\tau) \vert^2 \\
& \phantom{\leq}\; + 16 \Vert b \Vert^2 \sup_{\tau\in[0,t]} e^{-2\kappa(t-\tau)} \vert \Delta \dot{u}(\tau) \vert^2 
\end{align*}
for all $t \geq 0$. Since $M \geq N$, we have $\sum_{n \geq M+1} \vert \Delta v_n(\tau) \vert^2 \leq S_v(\tau)$. Hence, we obtain from (\ref{eq: stability residual dynamics - prel 4}) that
\begin{align*}
& \int_0^t e^{-(2\kappa+\delta)(t-\tau)} \sum_{n \geq M+1} \vert \Delta v_n(\tau) \vert^2 \,\mathrm{d}\tau \leq \\
& \dfrac{\gamma_9}{\delta} e^{-2\kappa t} \sup_{\tau\in[-h_M,0]} \Vert \Delta\Phi(\tau) \Vert^2 
+ \dfrac{\gamma_9}{\delta} \sup_{\tau\in[0,t]} e^{-2\kappa(t-\tau)} \vert \Delta u(\tau) \vert^2 
\end{align*}
The two latter inequalities imply the existence of a constant $\gamma_{11} > 0$ such that
\begin{align}
& \dfrac{1}{\gamma_{11}} \sum_{n \geq M+1} \vert \lambda_n \vert \vert \Delta w_n(t) \vert^2 \leq \label{eq: regulation - intermediate equation} \\
& e^{- 2 \kappa t} \sum_{n \geq M+1} \vert \lambda_n \vert \vert \Delta w_n(0) \vert^2 
+ e^{-2\kappa t} \sup_{\tau\in[-h_M,0]} \Vert \Delta\Phi(\tau) \Vert^2 \nonumber \\
& + \sup_{\tau\in[0,t]} e^{-2\kappa(t-\tau)} \vert \Delta u(\tau) \vert^2
+ \sup_{\tau\in[0,t]} e^{-2\kappa(t-\tau)} \vert \Delta \dot{u}(\tau) \vert^2 \nonumber
\end{align}
for all $t \geq 0$. Since $\Delta W(0) \in D(\mathcal{A}_{c,0})$ and $\vert \lambda_n \vert \leq \vert \lambda_n \vert^2$ for all $n \geq M+1$, we note that $\sum_{n \geq M+1} \vert \lambda_n \vert \vert \Delta w_n(0) \vert^2 \leq \Vert \mathcal{A}_{c,0} \Delta W(0) \Vert^2 \leq 2 \Vert \mathcal{A}_{c} \Delta X(0) \Vert^2 + 2 \Vert \mathcal{A}_{c}L \Vert^2 \vert \Delta u(0) \vert^2$
with $\Delta X(0) = \Delta \Phi(0)$ and $\vert \Delta u(0) \vert \leq \Vert K \Vert \Vert \Delta Y_a(0) \Vert + \vert \Delta p(0) \vert \leq \Vert K \Vert (\Vert \Delta \Phi(0) \Vert + \vert \Delta\zeta_0(0) \vert) + \vert \Delta p(0) \vert$. To conclude the proof, it is sufficient to study the two last terms of (\ref{eq: regulation - intermediate equation}). The estimation of the term involving $\Delta u$ immediately follows from (\ref{eq: main theorem - estimate trajectory}). Hence, only the term involving $\Delta \dot{u}$ needs to be investigated. From (\ref{eq: dynamics of deviations - truncated model}) and (\ref{eq: dynamics of deviations - control input}), we have, for all $t \geq 0$, $\vert \Delta \dot{u}(t) \vert
\leq \Vert K \Vert  \Vert \Delta\dot{Y}_a(t) \Vert  +  \vert \Delta\dot{p}(t) \vert$ with $\Vert \Delta\dot{Y}_a(t) \Vert \leq \Vert A_K - cI \Vert \Vert \Delta Y_a(t) \Vert + \vert c \vert \Vert \Delta Y_a(t-h(t)) \Vert + \Vert B_a \Vert \vert \Delta p(t) \vert + \vert \Delta r(t) \vert$. The claimed conclusion now follows from $\Vert \Delta Y_a(\tau) \Vert \leq \Vert \Delta X(\tau) \Vert + \vert \Delta \zeta(\tau) \vert$ for $\tau \geq -h_M$ and (\ref{eq: main theorem - estimate trajectory}).
\qed

\begin{rem}
In the context of Theorem~\ref{thm: regulation} dealing with classical solutions, the stability result stated by Theorem~\ref{thm: stability} can be strengthen as follows. First, it can be shown similarly to \cite[Eq.~42]{prieur2018feedback} that $\Vert f' \Vert^2 = - \cot(\theta) \vert f(1) \vert^2 + \dfrac{b+c}{a} \Vert f \Vert^2 - \dfrac{1}{a} \sum_{n \geq 0} \lambda_n \vert \langle f , e_n \rangle \vert^2$ for any $f \in D(\mathcal{A}_{c,0})$. Considering classical solutions, we can apply this identity to $\Delta W(t) \in D(\mathcal{A}_{c,0})$ where we note that estimates of $\Vert \Delta W(t) \Vert$ and $\vert [\Delta W(t)](1) \vert$ are provided by Theorem~\ref{thm: stability} and Theorem~\ref{thm: regulation}, respectively, while the series $\sum_{n \geq 0} \vert \lambda_n \vert \vert \Delta w_n(t) \vert^2$ has been evaluated in the proof of Theorem~\ref{thm: regulation}. Since $\Delta X(t) = \Delta W(t) + L \Delta u(t) \in D(\mathcal{A}_c) \subset H^2(0,1)$ with $\Vert [ L \Delta u(t) ]' \Vert = \frac{2}{\sqrt{3}} \vert \Delta u(t) \vert$, we infer that $\Vert [\Delta X(t)]' \Vert$, and hence $\Vert \Delta X(t) \Vert_{H^1(0,1)}$, is upper bounded by a term similar (i.e., with different constants $\overline{C}_i$) to the right-hand side of (\ref{eq: main theorem - estimate regulation}). If we further make the assumptions of Corollary~\ref{cor: regulation}, we obtain that $X(t)$ converges in $H^1(0,1)$ norm and hence, by the continuous embedding $H^1(0,1) \subset \mathcal{C}^0([0,1])$, in $L^\infty$ norm to $X_e$ when $t \rightarrow + \infty$.
\end{rem}

\section{Robustness with respect to delay mismatches}\label{sec: robustness}

In the previous sections, we have assumed the perfect knowledge of the state-delay $h$. This was used to build the dynamics of the integral component $\zeta$ given by (\ref{eq: zeta dynamics}). In this section, we discuss the robustness of the proposed control strategy with respect to delay mismatches. Assume that we dispose of an estimate $\hat{h}$ of the actual delay $h$ such that $\vert \hat{h} - h \vert \leq \delta$ for some constant $\delta > 0$. In this case, we replace the integral component $\zeta$, originally defined by (\ref{eq: zeta dynamics}), by the following dynamics:
\begin{subequations}\label{eq: robustness - zeta hat dynamics}
\begin{align}
\dot{\zeta}(t) & = \sum\limits_{n = 0}^N x_n(t) e_n(1) + c \{ \zeta(t-\hat{h}(t)) - \zeta(t) \} \\
& \phantom{=}\; + \alpha u(t) - r(t) , \\
\zeta(\tau) & = \zeta_0(\tau) , \quad \tau \in [-h_M,0]
\end{align}
\end{subequations}
Assuming that $\hat{h}$ satisfies the same assumptions as $h$, the only difference comparing to the previous developments occurs in the study of the truncated model. More precisely, the closed-loop truncated model, originally given by (\ref{eq: dynamics of deviations - truncated model}) and (\ref{eq: dynamics of deviations - IC trunated model}), becomes:
\begin{subequations}\label{eq: robustness - closed-loop truncated model}
\begin{align}
\Delta \dot{Y}_a(t) & = A_K \Delta Y_a(t) + c \{ \Delta Y_a(t-h(t)) - \Delta Y_a(t) \} \nonumber \\
& \phantom{=}\; + E_c \{ \Delta Y_a(t-\hat{h}(t)) - \Delta Y_a(t-h(t)) \} \nonumber \\
& \phantom{=}\; + B_a \Delta p(t) + \Delta \Gamma(t)  \\
\Delta Y_a(\tau) & = \Delta Y_{\Phi,a}(\tau) , \quad \tau\in[-h_M,0] 
\end{align}
\end{subequations}
with $E_c = \mathrm{diag}(0,\ldots,0,c) \in \mathbb{R}^{(N+2) \times (N+2)}$. Provided a suitable choice of the feeback gain $K$, the existence of a maximal delay mismatch $\delta > 0$ such that (\ref{eq: robustness - closed-loop truncated model}) is exponentially ISS with respect to the exogenous signals $\Delta p$ and $\Delta r$ follows from the following lemma.

\begin{lem}
Let $N \geq 1$, $0 < h_m < h_M$, $\mathbf{A},\mathbf{E} \in \R^{N \times N}$ and $c \in \R$. Assume that $\mathbf{A}$ is Hurwitz with simple eigenvalues $\mu_1,\ldots,\mu_N \in \C$ such that $\operatorname{Re}\mu_n < -3 \vert c \vert$ for all $1 \leq n \leq N$. Then there exist constants $\delta,\sigma,C_6,C_7 > 0$ such that, for any $x_0 \in \mathcal{C}^0([-h_M,0];\R^N)$, any $h_i \in \mathcal{C}^0(\R_+)$ with $i \in\{1,2,3\}$, $h_m \leq h_i(t) \leq h_M$, and $\vert h_2 - h_3 \vert \leq \delta$, and any $\tilde{q} \in \mathcal{C}^0(\R_+;\R^N)$, the trajectory of
\begin{subequations}\label{eq - robustness truncated model}
\begin{align}
\dot{x}(t) & = \mathbf{A} x(t) + c \left\{ x(t-h_1(t)) - x(t) \right\} \\
& \phantom{=}\; + \mathbf{E} \left\{ x(t-h_2(t)) - x(t-h_3(t)) \right\} + \tilde{q}(t) \nonumber , \\
x(\tau) & = x_0(\tau) , \quad \tau \in [-h_M,0]
\end{align}
\end{subequations}
satisfies the estimate $\Vert x(t) \Vert \leq C_6 e^{-\sigma t} \sup\limits_{\tau\in[-h_M,0]} \Vert x_0(\tau) \Vert + C_7 \sup\limits_{\tau\in[0,t]} e^{-\sigma (t-\tau)} \Vert \tilde{q}(\tau) \Vert$ for all $t \geq 0$.
\end{lem}

Hence, proceeding exactly as in the previous sections, we obtain the existence of a constant $\delta > 0$ such that, when replacing the definition (\ref{eq: zeta dynamics}) of the integral component $\zeta$ by (\ref{eq: robustness - zeta hat dynamics}), the conclusions of Theorems~\ref{thm: stability} and~\ref{thm: regulation} still hold true\footnote{With constants $\overline{C}_i$ of the estimates (\ref{eq: main theorem - estimate trajectory}) and (\ref{eq: main theorem - estimate regulation}) that are independent of a particularly selected $\hat{h}$.} for any estimated state-delay $\hat{h}$ satisfying the same assumptions as $h$ and with $\vert \hat{h} - h \vert \leq \delta$.

\textbf{Proof.}
Introducing $v_1(t) = x(t) - x(t-h_1(t))$, $v_2(t) = x(t-h_2(t)) - x(t-h_3(t))$, and $q(t) = \mathbf{E} v_2(t) + \tilde{q}(t)$, we obtain from Lemma~\ref{lem: prel lemma truncated model} that (\ref{eq: prel lemma truncated model - ISS}) holds. Since $\mathbf{A}$ is Hurwitz, we can assume that the constant $\sigma > 0$ involved in the latter equation (\ref{eq: prel lemma truncated model - ISS}) is further selected such that $\Vert e^{\mathbf{A}t} \Vert \leq M e^{-\sigma t}$ for all $t \geq 0$ and for some constant $M \geq 1$. Integrating (\ref{eq - robustness truncated model}) over either $[t-h_3(t),t-h_2(t)]$ if $h_3(t) \geq h_2(t)$ or $[t-h_2(t),t-h_3(t)]$ if $h_2(t) \geq h_3(t)$ for $t \geq h_M$, and combining these estimates, we obtain that 
\begin{align*}
\sup_{\tau \in [h_M,t]} e^{\sigma\tau} \Vert v_2(\tau) \Vert 
& \leq \{ e^{\delta \Vert \mathbf{A} \Vert} -1 \} e^{\sigma h_M} \sup_{\tau \in [0,t]} e^{\sigma\tau} \Vert x(\tau) \Vert \\
& \phantom{\leq}\; + \delta M \vert c \vert e^{\sigma h_M} \sup_{\tau\in[0,t]} e^{\sigma\tau} \Vert v_1(\tau) \Vert \\
& \phantom{\leq}\; + \delta M \Vert \mathbf{E} \Vert e^{\sigma h_M} \sup_{\tau\in[0,t]} e^{\sigma\tau} \Vert v_2(\tau) \Vert \\
& \phantom{\leq}\; + \delta M e^{\sigma h_M} \sup_{\tau\in[0,t]} e^{\sigma\tau} \Vert \tilde{q}(\tau) \Vert 
\end{align*}
for all $t \geq h_M$. From (\ref{eq: prel lemma truncated model - ISS}), the identity $v_1(t) = x(t) - x(t-h_1(t))$, and based on a small gain argument, we can fix $\delta > 0$ small enough (independently of $x_0$, $h_i$, and $\tilde{q}$) to obtain the existence of a constant $\gamma_{12} > 0$ such that $\sup_{\tau \in [h_M,t]} e^{\sigma\tau} \Vert v_2(\tau) \Vert \leq \gamma_{12} \sup_{\tau \in [-h_M,0]} \Vert x_0(\tau) \Vert + \gamma_{12} \sup_{\tau \in [0,h_M]} e^{\sigma\tau} \Vert v_2(\tau) \Vert + \gamma_{12} \sup_{\tau \in [0,t]} e^{\sigma\tau} \Vert \tilde{q}(\tau) \Vert$ for all $t \geq h_M$. Recalling that $v_2(t) = x(t-h_2(t)) - x(t-h_3(t))$, one can estimate for $t \in [0,h_M]$ the term $\sup_{\tau \in [0,t]} e^{\sigma\tau} \Vert v_2(\tau) \Vert$ from (\ref{eq - robustness truncated model}) and the use of Gr{\"o}nwall's inequality. Combining with the latter estimate, we obtain the existence of a constant $\gamma_{13} > 0$ such that $\sup_{\tau \in [0,t]} e^{\sigma\tau} \Vert v_2(\tau) \Vert \leq \gamma_{13} \sup_{\tau \in [-h_M,0]} \Vert x_0(\tau) \Vert + \gamma_{13} \sup_{\tau \in [0,t]} e^{\sigma\tau} \Vert \tilde{q}(\tau) \Vert$ for all $t \geq 0$. Since $q(t) = \mathbf{E} v_2(t) + \tilde{q}(t)$, the substitution of the latter estimate into (\ref{eq: prel lemma truncated model - ISS}) completes the proof.
\qed

\section{Simulation results}\label{sec: simulation}

We set $a = 0.2$, $b =2$, $c = 1$, and $\theta = \pi/3$. The first eigenvalues of $\mathcal{A}_{c,0}$ are approximately given by $\lambda_0 \approx 2.301$, $\lambda_1 \approx -1.668 > -2\sqrt{5}\vert c \vert$, and $\lambda_2 \approx -9.567 < -2\sqrt{5}\vert c \vert$. Hence we set $N = 1$. The feedback gain $K \in \R^{1 \times 3}$ is computed such that $A_K = A_a + B_a K$ is Hurwitz with simple eigenvalues $\mu_1 = -4$, $\mu_2 = -5$, $\mu_3 = -6$, selected so that $\mu_n < - 3 \vert c \vert$. The initial conditions of the plant and the integral component are set as $\phi(\tau,x) = 10 \cos(3 \pi \tau) x(1-x)^2$ and $\zeta_0(\tau) = \cos(3 \pi \tau) \zeta_a$ where $\zeta_a \in\R$ is selected such that (\ref{eq: classical solutions - compatibility condition}) holds. The numerical scheme consists of the modal approximation of the reaction–diffusion equation using its first 40 modes.

\begin{figure}
     \centering
     	\subfigure[State $y(t,x)$]{
		\includegraphics[width=3in]{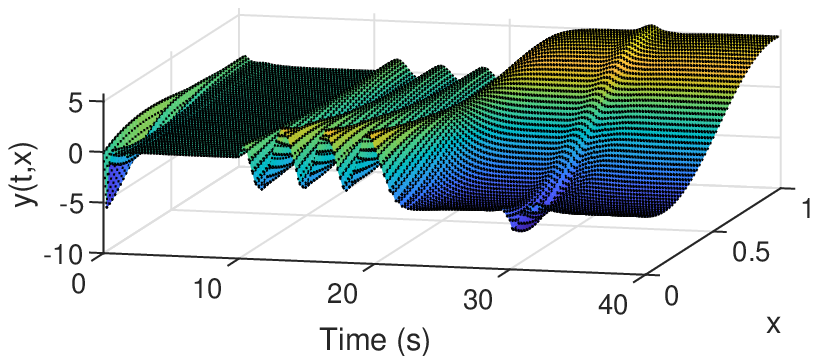}
		\label{fig: sim1 - state}
		}
     	\subfigure[Regulated output $z(t) = y(t,1)$]{
		\includegraphics[width=3in]{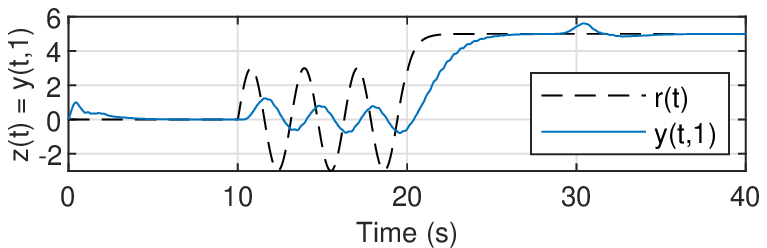}
		\label{fig: sim1 - output}
		}
     	\subfigure[Control input $y(t,0) = u(t)$]{
		\includegraphics[width=3in]{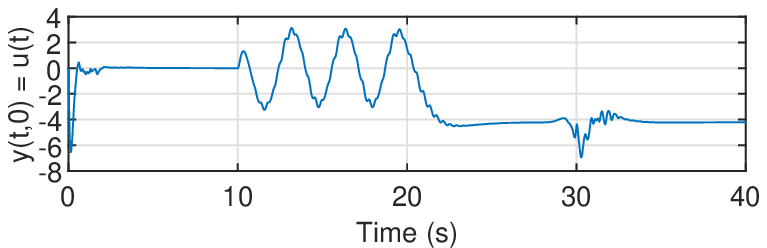}
		\label{fig: sim1 - input}
		}
     	\subfigure[Boundary perturbation $p(t)$]{
		\includegraphics[width=3in]{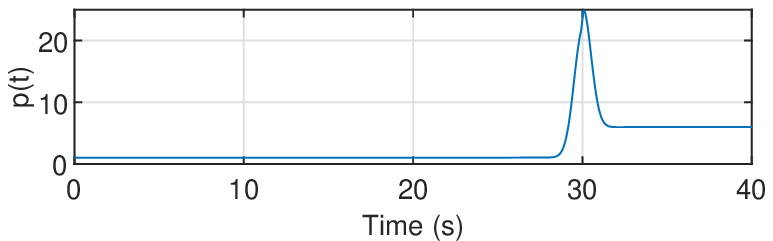}
		\label{fig: sim1 - perturbation}
		}
     \caption{Time evolution of the closed-loop system}
     \label{fig: sim1}
\end{figure}

The behavior of the closed-loop system composed of (\ref{eq: heat equation - ACS}), (\ref{eq: zeta dynamics}), and (\ref{eq: control input u}) is illustrated for the time varying delay $h(t) = 1 + \frac{1}{2}\sin(5\pi t+\pi/4)$ and the boundary perturbation $p(t)$ as shown in Fig.~\ref{fig: sim1 - perturbation}. The results are depicted in Fig.~\ref{fig: sim1}. During the 10 first seconds we observe that the control law achieves the stabilization of the closed-loop system: both the state and the regulated output converge to zero in spite of a constant perturbation $p(t) = 1$. Then, in order to evaluate the setpoint tracking capabilities of the system ouput (see Thm.~\ref{thm: regulation}), the reference signal is set as $r(t) = 5$ for $t > 20\,\mathrm{s}$ after an oscillatory transient. In conformity with the tracking estimate (\ref{eq: main theorem - estimate regulation}), we observe that the control strategy ensures the setpoint tracking of the reference signal $r(t)$ by the right Dirichlet trace $y(t,1)$. Around $t = 30\,\mathrm{s}$, the boundary perturbation $p(t)$ increases to reach (approximately) the value of 25 and then decreases to converge to the value of 6. It is seen that the impact of this perturbation on both the state trajectory and the regulated output are successfully eliminated due to the presence of the integral component. 

\begin{figure}
     \centering
     	\subfigure[Norm of the state $y(t,\cdot)$]{
		\includegraphics[width=3in]{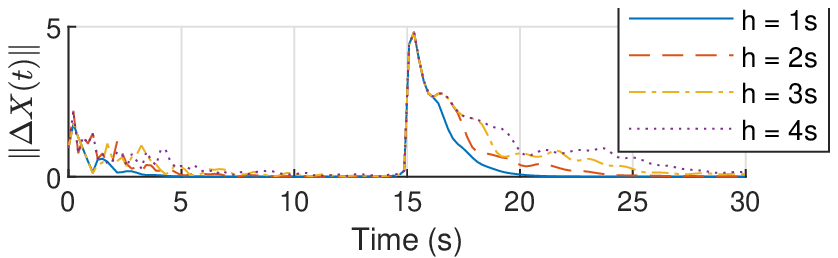}
		\label{fig: sim2 - state}
		}
     	\subfigure[Regulated output $z(t) = y(t,1)$]{
		\includegraphics[width=3in]{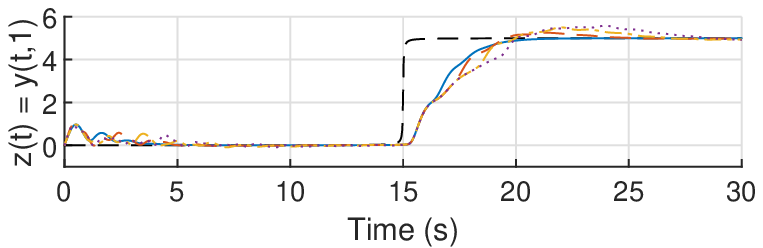}
		\label{fig: sim2 - output}
		}
     	\subfigure[Control input $y(t,0)= u(t)$]{
		\includegraphics[width=3in]{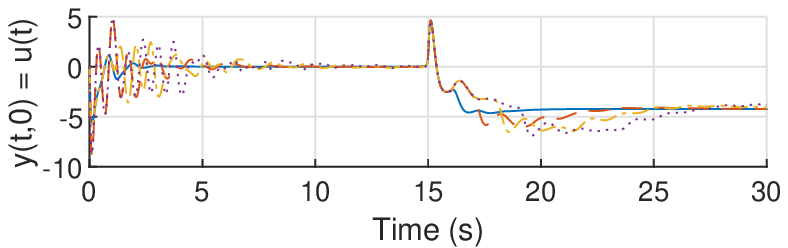}
		\label{fig: sim2 - input}
		}
     \caption{Impact of delay mismatches}
     \label{fig: sim2}
\end{figure}

Finally, Fig.~\ref{fig: sim2} illustrates the impact of delay mismatches on the closed-loop system performance. Here we set $\hat{h} = 1$ while considering increasing values for the actual delay $h \in\{1,2,3,4\}$. The boundary perturbation is set as $p = 1$. As expected, we observe a smooth degradation of the performances of the resulting closed-loop system.

\section{Conclusion}\label{sec: conclusion}

This paper has investigated the boundary PI regulation control of a reaction-diffusion equation in the presence of a state-delay in the reaction term. Our modal-based approach ensures the stability of the resulting closed-loop system as well as the setpoint regulation of the right Dirichlet trace. Future research directions may be concerned with extensions to the PI regulation control of either linear wave equations or semilinear heat equations in the presence of a state-delay.


\bibliographystyle{plain}        
\bibliography{autosam}           



\appendix

\end{document}